\documentclass[11pt]{article}
\usepackage{mathrsfs}
\usepackage{amsmath}
\usepackage{amsmath,amssymb,amsfonts,amsthm,fancyhdr}
\usepackage{epsfig,graphicx,picins,picinpar,subfigure}
\usepackage{pstricks}
\usepackage{fancyvrb}
\usepackage[numbers,sort&compress]{natbib}
\begin{document}

\title{\textbf{A new inequality involving primes}}
\author{\emph{Shaohua Zhang}$^{1,2}$}
\date{{\small 1 School of Mathematics, Shandong University,
Jinan,  Shandong, 250100, China \\
2 The key lab of cryptography technology and information security,
Ministry of Education, Shandong University, Jinan, Shandong, 250100,
China\\E-mail: shaohuazhang@mail.sdu.edu.cn}}
 \maketitle

\begin{abstract}
In this note, we find a new inequality involving primes and deduce
several Bonse-type inequalities.

\vspace{3mm}\textbf{Keywords:} prime, the $n$th prime, Bonse's
inequality, P\'{o}sa's inequality,  Panaitopol's inequality

\vspace{3mm}\textbf{2000 MR  Subject Classification:}\quad 11A41

\end{abstract}

\section{Introduction}
Denote the $n$th prime by $p_n$. In 1907, Bonse [1, 2] found and
proved two interesting inequalities which states that for $n\geq 4$,
$\prod_{i=1}^{i=n}p_i>p_{n+1}^2$ and for $n\geq 5$,
$\prod_{i=1}^{i=n}p_i>p_{n+1}^3$. Based on the first inequality, he
showed that a well known result which states that 30 is the largest
integer $N$ with the property that every integer $a$ with $1 < a <
N$ and $(a,N) = 1$ is prime. (This result has been further
generalized. See [3, 4]. ) In 2007, Betts [5] obtained  the
inequality $p_{k+1}-p_k< p_k(p_1p_2\cdots
p_{k-1}-p_k)/(p_{k+1}-p_k)$ by using Bonse's first inequality. Thus,
naturally, people are interesting in Bonse's inequalities. More
precisely, people would like to consider the inequalities about the
product of the first $n$ primes.

\vspace{3mm}In 1960, P\'{o}sa [6] refined firstly Bonse's
inequalities. He proved that for every integer $k>1$ there is an
$n_k$ such that $p_{n+1}^k<p_1p_2\cdots p_n$ for all $n>n_k$.
Moreover, the analogues of 30 are computed for the first few values
of $k$. In 1962, Mamangakis [7] proved that for $n\geq 11$,
$\prod_{i=1}^np_i>p_{4n}$ and for $n\geq 46$,
$\prod_{i=1}^{4n-9}p_i>p_{4n}^4$. In 1971, Reich [8] showed that for
every natural number $k$ there exists a natural number $N(k)$ such
that $\prod_{i=1}^np_i>p^2_{n+k}$ for all $n \geq N(k)$. In 1988,
S\'{a}ndor [9] showed that for $n\geq 3$, $ p_1p_2\cdots p_n\geq
p_1p_2\cdots p_{n-1}+p_n+p_{p_n-2}$, and for $n\geq 24$, $
p_1p_2\cdots p_n\geq p_{n+5}^2+p_{[n/2]}^2$, and for $n\geq 63$, $
p_1p_2\cdots p_n\geq p_{n+3}^3+p_{[n/3]}^6$, and so on. This refined
the Bonse's inequalities again. However his approach is quite
different from Bonse's. In 2000, using the Rosser-Schoenfeld and
Robin estimates, Panaitopol [10] proved that $p_1p_2\cdots
p_n>p_{n+1}^{n-\pi(n)}$, for all $n\geq 2$, where $\pi(n)$ is the
prime-counting function. In this note, we proved the following new
inequality involving primes:

\vspace{3mm}\noindent {\bf Theorem 1:~~}%
For integer $r\geq 20$, $p_{r+1}^{r-\pi(r)}>2^{p_{r+1}}$ and for
$1\leq r<20$, $p_{r+1}^{r-\pi(r)}<2^{p_{r+1}}$.

\vspace{3mm} By Theorem 1 and Panaitopol's inequality, we can deduce
the following result:

\vspace{3mm}\noindent {\bf Corollary 1:~~}%
For integer $r\geq 10$,  $p_1p_2\cdots p_r>2^{p_{r+1}}$. For integer
$0<r<10$ with $r\neq 8$, $p_1p_2\cdots p_r<2^{p_{r+1}}$.

\vspace{3mm} Corollary 1 improves P\'{o}sa's inequality in the
following form: for given integer $k\geq 5$, $p_1p_2\cdots
p_n>p^k_{n+1}$ for $n\geq 2k$. Bluntly speaking, the author likes
P\'{o}sa's inequality. In [11], using P\'{o}sa's inequality, the
author proved that there exists a prime $q $ such that for all prime
$p>q$ , if $1\leq a<p$, and $r$ is the smallest prime satisfying
$(r,a) = 1$, then $ 4r^3 < p$.

\vspace{3mm} Based on  Corollary 1, one could also get easily
several Bonse-type inequalities for the first few values of $n$. For
example, $\prod_{i=1}^{i=n}p_i>p_{n+1}^6$ provided $n\geq 10$, and
$\prod_{i=1}^{i=n}p_i>p_{n+1}^5$ provided $n\geq 8$.

\section{The Proof of Main Results}
\vspace{3mm}\noindent {\bf Lemma 1 [12]:~~}%
For $x>1$, $\pi(x)< \frac{1.25506 x}{\log x}$.

\vspace{3mm}\noindent {\bf Corollary 2:~~}%
For integer $r\geq 55$, $r-\pi(r)> (r+1)\log 2$.

\vspace{3mm}\noindent {\bf Proof:~~}%
Firstly, we can check directly that for $63\leq r\leq 149$,
$0.3r\geq \pi(r)+0.7$. If $r\geq 149$, then $\log r>5$ and
$7/r<0.05$. Therefore, $\frac{12.5506}{\log r}+\frac{7}{r}<3$. But,
by Lemma 1, $\pi(r)< \frac{1.25506 r}{\log r}$. So, $0.3r\geq
\pi(r)+0.7$, and for $r\geq 63$, $r-\pi(r)> (r+1)\times
0.7>(r+1)\log 2$. When $62\geq r\geq 55$, one can check directly
$r-\pi(r)> (r+1)\log 2$. This completes the proof of Corollary 2.

\vspace{3mm}\noindent {\bf Lemma 2 [12]:~~}%
For $x\geq 17$, $\pi(x)> \frac{x}{\log x}$.

\vspace{3mm}\noindent {\bf Proof of Theorem 1:~~}%
By Lemma 2, for  $p_{r+1}\geq 55$, $\pi(p_{r+1})>
\frac{p_{r+1}}{\log p_{r+1}}$. Hence, $r+1> \frac{p_{r+1}}{\log
p_{r+1}}$. By Corollary 2,  for integer $r\geq 55$, $r-\pi(r)>
(r+1)\log 2$.  Note also that for $r\geq 55$, $p_{r+1}\geq 55$. So,
$r-\pi(r)> \frac{p_{r+1}}{\log p_{r+1}}\log 2$, and for $r\geq 55$,
$p_{r+1}^{r-\pi(r)}>2^{p_{r+1}}$. When $20\leq r\leq 54$, one can
check directly that Theorem 1 is true  as follows:
$$ (20-\pi (20))\log p_{21}=12\times \log 73>51.4>73 \times 0.7> 73\log 2$$
$$ (21-\pi (21))\log p_{22}=13\times \log 79>56.8>79 \times 0.7> 79\log 2$$
$$ (22-\pi (22))\log p_{23}=14\times \log 83>61.8>83 \times 0.7> 83\log 2$$
$$ (23-\pi (23))\log p_{24}=14\times \log 89>62.8>89 \times 0.7> 89\log 2$$
$$ (24-\pi (24))\log p_{25}=15\times \log 97>68.6>97 \times 0.7> 97\log 2$$
$$ (25-\pi (25))\log p_{26}=16\times \log 101>73.8>101 \times 0.7>101\log 2$$
$$ (26-\pi (26))\log p_{27}=17\times \log 103>78.7>103 \times 0.7>103\log 2$$
$$ (27-\pi (27))\log p_{28}=18\times \log 107>84.1>107 \times 0.7>107\log 2$$
$$ (28-\pi (28))\log p_{29}=19\times \log 109>89.1>109 \times 0.7>109\log 2$$
$$ (29-\pi (29))\log p_{30}=19\times \log 113>89.8>113 \times 0.7>113\log 2$$
$$ (30-\pi (30))\log p_{31}=20\times \log 127>96.8>127 \times 0.7>127\log 2$$
$$ (31-\pi (31))\log p_{32}=20\times \log 131>97.5>131 \times 0.7>131\log 2$$
$$ (32-\pi (32))\log p_{33}=21\times \log 137>103.3>137 \times 0.7>137\log 2$$
$$ (33-\pi (33))\log p_{34}=22\times \log 139>108.5>139 \times 0.7>139\log 2$$
$$ (34-\pi (34))\log p_{35}=23\times \log 149>115.0>149 \times 0.7>149\log 2$$
$$ (35-\pi (35))\log p_{36}=24\times \log 151>120.4>151 \times 0.7>151\log 2$$
$$ (36-\pi (36))\log p_{37}=25\times \log 157>126.4>157 \times 0.7>157\log 2$$
$$ (37-\pi (37))\log p_{38}=25\times \log 163>127.3>163 \times 0.7>163\log 2$$
$$ (38-\pi (38))\log p_{39}=26\times \log 167>133.0>167 \times 0.7>167\log 2$$
$$ (39-\pi (39))\log p_{40}=27\times \log 173>139.1>173 \times 0.7>173\log 2$$
$$ (40-\pi (40))\log p_{41}=28\times \log 179>145.2>179 \times 0.7>179\log 2$$
$$ (41-\pi (41))\log p_{42}=28\times \log 181>145.5>181 \times 0.7>181\log 2$$
$$ (42-\pi (42))\log p_{43}=29\times \log 191>152.3>191 \times 0.7>191\log 2$$
$$ (43-\pi (43))\log p_{44}=30\times \log 193>152.6>193 \times 0.7>193\log 2$$
$$ (44-\pi (44))\log p_{45}=30\times \log 197>158.4>197 \times 0.7>197\log 2$$
$$ (45-\pi (45))\log p_{46}=31\times \log 199>164.0>199 \times 0.7>199\log 2$$
$$ (46-\pi (46))\log p_{47}=32\times \log 211>171.2>211 \times 0.7>211\log 2$$
$$ (47-\pi (47))\log p_{48}=32\times \log 223>173.0>223 \times 0.7>223\log 2$$
$$ (48-\pi (48))\log p_{49}=33\times \log 227>179.0>227 \times 0.7>227\log 2$$
$$ (49-\pi (49))\log p_{50}=34\times \log 229>184.7>229 \times 0.7>229\log 2$$
$$ (50-\pi (50))\log p_{51}=35\times \log 233>190.7>233 \times 0.7>233\log 2$$
$$ (51-\pi (51))\log p_{52}=36\times \log 239>197.1>239 \times 0.7>239\log 2$$
$$ (52-\pi (52))\log p_{53}=37\times \log 241>202.9>241 \times 0.7>241\log 2$$
$$ (53-\pi (53))\log p_{54}=37\times \log 251>204.4>251 \times 0.7>251\log 2$$
$$ (54-\pi (54))\log p_{55}=38\times \log 257>210.8>257 \times 0.7>257\log 2$$

When $1\leq r<20$, one can
 check similarly $p_{r+1}^{r-\pi(r)}<2^{p_{r+1}}$. This completes the
proof of Theorem 1.

\vspace{3mm}\noindent {\bf Proof of Corollary 1:~~}%
By Panaitopol's inequality, for $r\geq 20$,  we have $p_1p_2\cdots
p_r>2^{p_{r+1}}$. When $10\leq r\leq 19$, one can check directly
that Corollary 1 is true. The remaining case can be checked
similarly. This completes the proof of  Corollary 1.

\vspace{3mm} Finally, we prove that Corollary 1 improves P\'{o}sa's
inequality in the following form: for given integer $k\geq 5$,
$p_1p_2\cdots p_n>p^k_{n+1}$ for $n\geq 2k$. Note that for $k\geq
5$, we have $n\geq 2k\geq 10$. So, by Corollary 1, we have
$p_1p_2\cdots p_n>2^{p_{n+1}}$. But by Lemma 1, we have  $n+1<
\frac{1.25506 p_{n+1}}{\log p_{n+1}}$. On the other hand,
$\frac{n+1}{1.25506}>\frac{n}{2\log2}$ since
$\frac{1}{1.25506}>\frac{1}{2\log2}$. Thus, $\frac{p_{n+1}}{\log
p_{n+1}}>\frac{n}{2\log 2}$. So, $2^{p_{n+1}}>p^{n/2}_{n+1}$ and
$p_1p_2\cdots p_n>2^{p_{n+1}}>p^{n/2}_{n+1}\geq p^k_{n+1}$. This
completes the proof. As an application, one can deduce easily that
$\prod_{i=1}^{i=n}p_i>p_{n+1}^6$ provided $n\geq 10$, and
$\prod_{i=1}^{i=n}p_i>p_{n+1}^5$ provided $n\geq 8$.

\section{Acknowledgements}
Thank my advisor Professor Xiaoyun Wang for her valuable help. Thank
Institute for Advanced Study in Tsinghua University for providing me
with excellent conditions. This work was partially supported by the
National Basic Research Program (973) of China (No. 2007CB807902)
and the Natural Science Foundation of Shandong Province (No.
Y2008G23).

\clearpage
\end{document}